\documentclass{amsart}
\usepackage{amsthm}
\usepackage{amssymb}
\usepackage{upref}
\usepackage{amscd}

\numberwithin{equation}{section}

\newtheorem{theorem}[equation]{Theorem}
\newtheorem{lemma}[equation]{Lemma}
\newtheorem{proposition}[equation]{Proposition}
\newtheorem{corollary}[equation]{Corollary}

\theoremstyle{definition}

\DeclareMathOperator{\ann}{ann}

\DeclareMathOperator{\Hom}{Hom}

\DeclareMathOperator{\Ext}{Ext}

\DeclareMathOperator{\Tor}{Tor}

\DeclareMathOperator{\colim}{colim}

\DeclareMathOperator{\im}{im}

\newcommand{\cat}[1]{\mathcal{#1}}

\newcommand{\mathcolon}{\colon\,}

\newcommand{\ulp}{\textup{(}}
\newcommand{\urp}{\textup{)}}

\begin{document}
 
\title{The generating hypothesis in the derived category of a ring} 

\date{\today}

\author{Mark Hovey}
\address{Department of Mathematics \\ Wesleyan University
\\ Middletown, CT 06459}
\email{hovey@member.ams.org}

\author{Keir Lockridge}
\address{Department of Mathematics \\ University of Washington
\\ Box 354350 \\
Seattle, WA 98195}
\email{lockridg@math.washington.edu}

\author{Gena Puninski}
\address{Department of Mathematics \\
University of Manchester \\
Booth Street East \\
Manchester M13 PL \\
United Kingdom}
\email{gpuninski@maths.man.ac.uk}


\begin{abstract}
We show that a strong form (the fully faithful version) of the
generating hypothesis, introduced by Freyd in algebraic topology, holds
in the derived category of a ring $R$ if and only if $R$ is von Neumann
regular.  This extends results of the second author~\cite{lockridge}.
We also characterize rings for which the original form (the faithful
version) of the generating hypothesis holds in the derived category of
$R$.  These must be close to von Neumann regular in a precise sense,
and, given any of a number of finiteness hypotheses, must be von Neumann
regular.  However, we construct an example of such a ring that is not
von Neumann regular, and therefore does not satisfy the strong form of
the generating hypothesis.
\end{abstract}

\maketitle

\section*{Introduction}

The generating hypothesis was introduced by Peter
Freyd~\cite{freyd-stable} in algebraic topology, where it is the
assertion that any map $f\mathcolon X\xrightarrow{}Y$ of finite spectra
that is $0$ on stable homotopy groups is in fact null homotopic.  The
generating hypothesis is widely considered to be one of the most
important and difficult problems in stable homotopy theory.  It has many
implications for the structure of the stable homotopy ring
$\pi_{*}S^{0}$ of the sphere, implying for example that it is totally
non-coherent~\cite{freyd-stable} and that the $p$-completion
$\pi_{*}S^{0}_{p}$ is a self-injective ring~\cite{hovey-generating}.
Somewhat surprisingly, Freyd proved that the generating hypothesis in
fact implies that the map
\[
[X,Y] \xrightarrow{} \Hom_{\pi_{*}S^{0}}(\pi_{*}X, \pi_{*}Y)
\]
from maps of finite spectra to maps of their stable homotopy modules is
not only injective but also surjective.  That is, the generating
hypothesis implies that the stable homotopy functor is fully faithful on
finite spectra.  

One approach to understanding the generating hypothesis is to look at
analogous questions in other categories.  Following the second
author~\cite{lockridge}, we say that a ring $R$ \textbf{satisfies the
generating hypothesis} if whenever $f\mathcolon X\xrightarrow{}Y$ is a
map of perfect complexes in the derived category $\cat{D}(R)$ of $R$ and
$H_{*}f=0$, then $f=0$.  Recall that a perfect complex is a bounded
chain complex of finitely generated projective (right) modules, and that
$f=0$ in $\cat{D}(R)$ exactly when $f$ is chain homotopic to $0$ (for
maps of perfect complexes).  Perfect complexes are the algebraic
analogue of finite spectra, as they are the small objects in
$\cat{D}(R)$.  Thus $R$ satisfies the generating hypothesis exactly when
the homology functor is faithful on perfect complexes.  Let us also say
that $R$ \textbf{satisfies the strong generating hypothesis} if the
homology functor is fully faithful on perfect complexes.

The second author noticed~\cite[Section~4]{lockridge} that the homology
functor is faithful on all of $\cat{D}(R)$ if and only if all right
$R$-modules are projective; that is, if and only if $R$ is semisimple.
Since perfect complexes are the small objects of $\cat{D}(R)$ and
finitely presented modules are the small $R$-modules, it is natural to
conjecture (as the second author did in~\cite{lockridge}) that the
homology functor is faithful on perfect complexes (that is, $R$
satisfies the generating hypothesis) if and only if all finitely
presented right $R$-modules are projective; that is, if and only if $R$
is von Neumann regular.  The second author verified that all von Neumann
regular rings do satisfy the generating hypothesis, and proved that if
$R$ satisfies the generating hypothesis and is either commutative or
right coherent, then $R$ is von Neumann regular~\cite{lockridge}.

In this paper, we first prove that $R$ satisfies the strong generating
hypothesis if and only if $R$ is von Neumann regular.  We then consider
the generating hypothesis, in effect asking whether the generating
hypothesis implies the strong generating hypothesis.  We prove that $R$
satisfies the generating hypothesis if and only if all short exact
sequences of finitely presented modules split, and all submodules of
flat modules are flat.  This makes $R$ close to von Neumann regular, and
in fact if $R$ is local or satisfies one of several finiteness
hypotheses it forces $R$ to be von Neumann regular.  However, we 
construct an example of a ring that satisfies the generating hypothesis
but is not von Neumann regular.  Over this ring, then, the homology
functor is faithful on perfect complexes but not full.  

The authors would like to thank Grigory Garkusha for many helpful
discussions.  

All $R$-modules $M$ will be right $R$-modules in this paper, so that,
for example, $\cat{D}(R)$ is the unbounded derived category of right
$R$-modules.  The differential $d$ in a chain complex $P$ will lower
dimension, so that $d_{n}\mathcolon P_{n}\xrightarrow{}P_{n-1}$.  We
will denote $\ker d_{n}$ by $Z_{n}P$ and $\im d_{n}$ by $B_{n-1}P$. If
$M$ is an $R$-module, then $D^{n}(M)$ denotes the complex which is $M$
in degree $n$ and $n-1$ and $0$ elsewhere, with $d_{n}$ being the
identity.  $S^{n}M$ denotes the complex that is $M$ in degree $n$ and $0$
elsewhere.  

\section{The strong generating hypothesis}\label{sec-strong}

We begin by recalling the second author's characterization of semisimple
rings.  

\begin{lemma}\label{lem-split}
Suppose $P$ is a perfect complex of $R$-modules with both $B_{n}P$ and
$H_{n}P$ projective for all $n$.  Then $P\cong \bigoplus_{n}
S^{n}(H_{n}P)$ in $\cat{D}(R)$.  In this case, the natural map
$[P,Q]\xrightarrow{}\Hom_{R}(H_{*}P, H_{*}Q)$ is an isomorphism for all
complexes $Q$.  
\end{lemma}

\begin{proof}
We have 
\[
P_{n} \cong Z_{n}P \oplus B_{n-1}P \cong B_{n}P\oplus H_{n}P \oplus
B_{n-1}P.  
\]
>From this it follows that $P\cong \bigoplus_{n} D^{n}(B_{n-1}P) \oplus
S^{n}(H_{n}P)$, which is isomorphic to $\bigoplus_{n}S^{n}(H_{n}P)$ in
$\cat{D}(R)$.  

A chain map from $S^{n}H_{n}P$ to a complex $Q$ is the same thing as a
map $f\mathcolon H_{n}P\xrightarrow{}Z_{n}Q$, and such a map is chain
homotopic to $0$ exactly when there is a map $D\mathcolon
H_{n}P\xrightarrow{}Q_{n+1}$ such that $dD=f$.  Since $H_{n}P$ is
projective, $f$ is chain homotopic to $0$ if and only if $f$ lands in
$B_{n}Q$.  Using projectivity of $H_{n}P$ again, we conclude that
$[S^{n}H_{n}P,Q]\cong \Hom_{R}(H_{n}P,H_{n}Q)$.   
\end{proof}

\begin{proposition}\label{prop-semisimple}
A ring $R$ is semisimple if and only if the homology functor is faithful
on $\cat{D}(R)$.  Furthermore, in this case, the homology functor is in
fact fully faithful on $\cat{D}(R)$.  
\end{proposition}

\begin{proof}
Suppose the homology functor is faithful in $\cat{D}(R)$.  Take two
$R$-modules $M$ and $N$, and take a projective resolution $P_{*}$ of
$M$.  Then an element of $\Ext^{s}(M,N)$ is represented by a map from
$P_{*}$ to $N$, thought as a complex concentrated in degree $s$.  This
map is necessarily $0$ in homology when $s>0$.  Thus $\Ext^{s}(M,N)=0$
for all $s>0$ and all $M,N$, so every $R$-module is projective and $R$
is semisimple. 

On the other hand, if $R$ is semisimple, then Lemma~\ref{lem-split}
implies that homology is fully faithful.
\end{proof}

The analogue for the generating hypothesis is the following theorem.  

\begin{theorem}\label{thm-regular}
A ring $R$ satisfies the strong generating hypothesis if and only if $R$
is von Neumann regular.  In this case, the natural map 
\[
[P,Q] \xrightarrow{} \Hom_{R}(H_{*}P,H_{*}Q) 
\]
is an isomorphism for all perfect complexes $P$ and arbitrary complexes
$Q$.  
\end{theorem}

Recall that $R$ is von Neumann regular if and only if, for every $x\in
R$, there is a $y\in R$ with $x=xyx$.  The standard reference for von
Neumann regular rings is~\cite{goodearl}; the book~\cite{lam} takes an
approach based on module categories, so contains some different and
useful results about von Neumann regular rings.  A standard
characterization is that $R$ is von Neumann regular if and only if all
$R$-modules are flat, which is true if and only if all finitely
presented $R$-modules are projective.

\begin{proof}
Suppose $R$ satisfies the strong generating hypothesis.  Then,
\[
\ann_{\ell} \ann_{r} (Rx) = Rx
\]
for all $x\in R$, by~\cite[Proposition~2.7]{lockridge}.  Now take $x\in
R$, and consider the perfect complex $P$ with $P_{i}=R$ if $i=0,1$ and
$P_{i}=0$ otherwise, with the differential $P_{1}\xrightarrow{}P_{0}$
being left multiplication by $x$.  This complex has $H_{0}(P)=R/xR$ and
$H_{1}(P)=\ann_{r} (x)$.  By the strong generating hypothesis, there
exists a chain map $\phi \mathcolon P\xrightarrow{}P$ such that
$H_{1}(\phi)=0$ and $H_{0}(\phi)=1$, the identity of $R/xR$.
Translating, this means there exist elements $a,b\in R$ such that
$xa=bx$ with $a\in \ann_{\ell}\ann_{r}(x)$ (so that $H_{1}(\phi)=0$) and
$b=1+xc$ for some $c\in R$ (so that $H_{0}(\phi )=1$).  But then $a=dx$
for some $d\in R$, so we have
\[
xdx = xa = bx = (1+xc)x = x + xcx.  
\]
This means that $x=x(d-c)x$.  Since $x$ was arbitrary, $R$ is von
Neumann regular.  

Conversely, suppose $R$ is von Neumann regular, and $P$ is a perfect
complex.  In a von Neumann regular ring, finitely generated submodules
of projectives are projective~\cite[p.44]{lam}, so $B_{n}P$ is finitely
generated projective for all $n$.  Then $Z_{n}P$, as the kernel of the
(necessarily split) surjection $P_{n}\xrightarrow{}B_{n-1}P$, is also
finitely generated projective for all $n$.  Hence $H_{n}P$ is finitely
presented, and so is projective for all $n$.  Now Lemma~\ref{lem-split}
implies that homology is fully faithful on maps out of perfect
complexes.
\end{proof}

Recall from~\cite{lockridge} that if $R$ is either commutative or right
coherent and $R$ satisifes the generating hypothesis, then $R$ is von
Neumann regular.  Hence we get the following corollary.  

\begin{corollary}\label{cor-regular}
If $R$ is either commutative or coherent, then $R$ satisfies the
generating hypothesis if and only if $R$ satisfies the strong
generating hypothesis.  
\end{corollary}

The second author also investigated the generating hypothesis from the
viewpoint of global stable homotopy theory.  Using the results
of~\cite{lockridge}, we get the following corollary.

\begin{corollary}\label{cor-thick}
A ring $R$ satisfies the strong generating hypothesis if and only if, in
$\cat{D}(R)$, the thick subcategory generated by $R$ is the collection
of retracts of finite coproducts of suspensions of $R$.  
\end{corollary}

Recall that a full subcategory of a triangulated category is called
\textbf{thick} if it is closed under shifts, retracts, and cofibers; the
thick subcategory generated by $R$ consists of the perfect complexes.
This corollary follows from~\cite[Proposition~5.1]{lockridge}, and
indicates how different stable homotopy theory must be from the derived
category of a ring if the generating hypothesis in stable homotopy is to
be true, since there are many finite spectra that are not retracts of
finite coproducts of suspensions of the sphere.

\section{Rings that satisfy the generating hypothesis}\label{sec-class}

Having dealt with the strong generating hypothesis, we now turn our
attention to the generating hypothesis.  The object of this section to
prove the following theorem.  

\begin{theorem}\label{thm-class}
A ring $R$ satisfies the generating hypothesis if and only if $R$ has
weak global dimension at most $1$ and all finitely presented $R$-modules
are FP-injective.  
\end{theorem}

Weak global dimension at most $1$ is of course equivalent to the
statement that submodules of flat modules are flat.  Recall that a
module $M$ is said to be \textbf{FP-injective} if $\Ext^{1}(F,M)=0$ for
all finitely presented modules $F$; thus all finitely presented modules
are FP-injective if and only if all short exact sequences of finitely
presented modules split.  FP-injective modules seem to have been
introduced in~\cite{stenstrom-fp-injective}; a good guide to the
literature can be found in~\cite[Chapter~6]{faith-rings}.  An
FP-injective module is sometimes called \textbf{absolutely pure},
because $M$ is FP-injective if and only if every short exact sequence
\[
0 \xrightarrow{} M \xrightarrow{} N \xrightarrow{} P \xrightarrow{} 0
\]
is pure (that is, remains exact upon tensoring with any left
$R$-module).  See~\cite[Theorem~4.89(5)]{lam} for a proof of this
equivalence.  

To compare the rings of Theorem~\ref{thm-class} with von Neumann
regular rings, the following lemma is helpful.  

\begin{lemma}\label{lem-fp-inj}
A ring $R$ is von Neumann regular if and only if every $R$-module is
FP-injective.  
\end{lemma}

This lemma is well-known, but does not appear in~\cite{goodearl}
or~\cite{lam}, so we include the proof for the convenience of the
reader.  

\begin{proof}
Suppose $R$ is von Neumann regular, and $M$ is an $R$-module.  Choose a
short exact sequence $\cat{E}$
\[
0 \xrightarrow{} M \xrightarrow{} I \xrightarrow{} N \xrightarrow{} 0
\]
where $I$ is injective.  Since $N$ is necesarily flat, this sequence is
pure~\cite[Theorem~4.85]{lam}.  Hence, if $F$ is finitely presented,
$\Hom (F, \cat{E})$ is still exact~\cite[Theorem~4.89(5)]{lam}, and so
$\Ext^{1}(F,M)=0$ and $M$ is FP-injective.  

Conversely, if every module is FP-injective, another application
of~\cite[Theorem~4.89(5)]{lam} shows that every short exact sequence
is pure.  Then~\cite[Theorem~4.85]{lam} shows that every module is
flat, as required. 
\end{proof}

We now begin the proof of Theorem~\ref{thm-class}.  Our first task is to
characterize the homology groups of perfect complexes.  

\begin{proposition}\label{prop-small}
Suppose $R$ is a ring.  An $R$-module $M$ is a homology module of a
perfect complex of $R$-modules if and only if there exists a finitely
presented module $F$ such that $M$ embeds in $F$ and the quotient $F/M$
embeds in a projective module. Furthermore, in this case, there is a
perfect complex $P$ such that $P_{n}=0$ unless $n=0,1,2$ and
$M=H_{1}P$.  
\end{proposition}

\begin{proof}
Suppose $M=H_{n}P$, where each $P_{i}$ is a finitely generated
projective module.  Then we have a short exact sequence
\[
0 \xrightarrow{} M \xrightarrow{} P_{n}/B_{n}P \xrightarrow{d_{n}}
B_{n-1}P \xrightarrow{}0,
\]
$P_{n}/B_{n}P$ is finitely presented and $B_{n-1}P$ embeds in the
projective module $P_{n-1}$.

Conversely, suppose $M$ embeds in the finitely presented module $F$ and
the quotient $F/M$ embeds in the projective module $P_{0}$, which we can
assume is finitely generated since $F$ is so.  Choose a presentation 
\[
P_{2}\xrightarrow{d_{2}} P_{1}\xrightarrow{p} F \xrightarrow{} 0
\]
of $F$, where $P_{1}$ and $P_{2}$ are finitely generated projectives.
Define the map $d_{1}\mathcolon P_{1}\xrightarrow{}P_{0}$ to be the
composite 
\[
P_{1} \xrightarrow{p} F \xrightarrow{} F/M \xrightarrow{} P_{0}.
\]
This defines a three-term perfect chain complex $P$.  Pulling back
the presentation of $F$ through the inclusion $M\xrightarrow{}F$ shows
that $H_{1}P\cong M$.   
\end{proof}

\begin{corollary}\label{cor-small}
Suppose $P$ is a perfect complex.  Then each cycle module $Z_{n}P$ is a
homology module of some perfect complex.
\end{corollary}

\begin{proof}
Note that $Z_{n}P$ is a submodule of the finitely presented module
$P_{n}$ and the quotient $P_{n}/Z_{n}P$ embeds in the projective module
$P_{n-1}$.  
\end{proof}

We now take a significant step towards Theorem~\ref{thm-class} by
showing how FP-injective modules get involved.  

\begin{theorem}\label{thm-generating}
Let $R$ be a ring, and let $Q$ be an arbitrary object of $\cat{D}(R)$.
Then the generating hypothesis with target $Q$ is true in
$\cat{D}(R)$ if and only if $H_{n}Q$ is FP-injective for all $n$.
In particular, $R$ satisfies the generating hypothesis if and only if
all homology modules of perfect complexes are FP-injective.  
\end{theorem}

The \textbf{generating hypothesis with target $Q$} is the statement that
any map $f\mathcolon P\xrightarrow{}Q$ in $\cat{D}(R)$ where $P$ is a
perfect complex and $H_{*}f=0$ has $f=0$.  So $R$ satisfies the
generating hypothesis if and only if $R$ satisfies the generating
hypothesis with target $Q$ for all perfect complexes $Q$.

Note in particular that this theorem and Lemma~\ref{lem-fp-inj} imply
that $R$ satisfies the generating hypothesis with target $Q$ for all
(not necessarily perfect) $Q$, if and only if $R$ is von Neumann
regular.

\begin{proof}
Suppose first that the generating hypothesis with target $Q$ holds, and
consider a finitely presented module $F$ and an integer $n$.  Choose a
finite presentation
\[
P_{n}\xrightarrow{d_{n}} P_{n-1} \xrightarrow{} F \xrightarrow{} 0
\]
of $F$, so that, by letting $P_{i}=0$ for $i\neq n,n-1$, we get a
perfect complex $P_{*}$ with $H_{n-1}P_{*}=F$.  To prove that
\[
\Ext^{1}(F, H_{n}Q_{*})=0,
\]
it suffices to show that any map
\[
f\mathcolon P_{n}/Z_{n}P \xrightarrow{}H_{n}Q_{*}
\]
extends to a map $g\mathcolon P_{n-1}\xrightarrow{}H_{n}Q_{*}$ with
$g\overline{d_{n}}=f$, where $\overline{d_{n}}$ is the map induced by
$d_{n}$.

Since $P_{n}$ is projective, there is a map $\phi_{n}\mathcolon
P_{n}\xrightarrow{}Q_{n}$ such that the composite 
\[
P_{n}\xrightarrow{\phi_{n}} Q_{n} \xrightarrow{q} Q_{n}/B_{n}Q
\]
is the composite 
\[
P_{n} \xrightarrow{p} P_{n}/Z_{n}P \xrightarrow{f} H_{n}Q_{*}
\xrightarrow{i} Q_{n}/B_{n}Q.  
\]
Now let $\phi_{n-1}\mathcolon P_{n-1}\xrightarrow{}Q_{n-1}$ be the zero map.
Then $\phi \mathcolon P_{*}\xrightarrow{}Q_{*}$ is a chain map.  Indeed,
write $d_{n}\mathcolon Q_{n}\xrightarrow{}Q_{n-1}$ as $d_{n}=rq$.  Then 
\[
d_{n}\phi_{n} = rq\phi_{n} = rifp =0
\]
since $ri=0$.  Furthermore, $\phi$ induces the zero map on homology,
because if $x\in Z_{n}P$, then $q\phi_{n}x=0$, so $\phi_{n}x$ is a
boundary.

If the generating hypothesis is true, then $\phi$ must be chain
homotopic to $0$.  This gives us maps $D_{n-1}\mathcolon
P_{n-1}\xrightarrow{}Q_{n}$ and $D_{n}\mathcolon P_{n}\xrightarrow{}Q_{n+1}$
such that $d_{n}D_{n-1}=0$ and $D_{n-1}d_{n}+d_{n+1}D_{n}=\phi_{n}$.
Since $d_{n}\mathcolon Q_{n}\xrightarrow{}Q_{n-1}$ factors through
$Q_{n}/B_{n}Q$ as $d_{n}=\overline{d_{n}}q$, we conclude that
$\overline{d_{n}}qD_{n-1}=0$, so there exists a map $g\mathcolon
P_{n-1}\xrightarrow{}H_{n}Q_{*}$ such that $ig=qD_{n-1}$.  Of course, we
claim that $g$ is the desired extension.  To see this, apply $q$ to the
relation 
\[
D_{n-1}d_{n}+d_{n+1}D_{n}=\phi_{n}
\]
to get 
\[
qD_{n-1}d_{n} = ifp \text{ or } igd_{n} = ifp.
\]
Writing $d_{n}=\overline{d_{n}}p$ and using the fact that $i$ is a
monomorphism and $p$ is an epimorphism, we conclude that
$g\overline{d_{n}}=f$, as required.  

Now suppose that every homology group of $Q$ is FP-injective, and $\phi
\mathcolon P_{*}\xrightarrow{}Q_{*}$ is a map of perfect complexes that
induces $0$ on homology.  We will construct a chain homotopy
$D_{n}\mathcolon P_{n}\xrightarrow{}Q_{n+1}$ such that
$d_{n+1}D_{n}+D_{n-1}d_{n}=\phi_{n}$ by induction on $n$.  Our induction
hypothesis will be that we have constructed $D_{i}$ for $i\leq n-1$ and
that $\phi_{n}-D_{n-1}d_{n}$, which is a map from $P_{n}$ to $Q_{n}$, in
fact lands in the boundaries $B_{n}Q$.  Getting started is easy since
$P$ is bounded below.  For the induction step, our hypothesis gives us
the commutative square below,
\[
\begin{CD}
P_{n+1}/Z_{n+1}P @>\overline{\phi_{n+1}}>> Q_{n+1}/B_{n+1}Q \\
@Vd_{n+1}VV @VVd_{n+1}V \\
P_{n} @>>\phi_{n} - D_{n-1}d_{n}> B_{n}Q
\end{CD}
\]
where $\overline{\phi_{n+1}}$ exists because $\phi$ is zero on homology,
so must take cycles to boundaries.  We will construct a lifting
$\overline{D_{n}}\mathcolon P_{n}\xrightarrow{}Q_{n+1}/B_{n+1}Q$ in this
square.  First of all, there is obviously a map
\[
E_{n}\mathcolon P_{n}\xrightarrow{}Q_{n+1}/B_{n+1}Q
\]
such that $d_{n+1}E_{n}=\phi_{n}-D_{n-1}d_{n}$, simply because $P_{n}$
is projective.  Then
\[
d_{n+1}(\overline{\phi_{n+1}}-E_{n}d_{n+1}) = d_{n+1}\overline{\phi_{n+1}}
-\phi_{n}d_{n+1}+D_{n-1}d_{n}d_{n+1} =0.
\]
Hence $\overline{\phi_{n+1}}-E_{n}d_{n+1}$ is a map from $P_{n+1}/Z_{n+1}P$
to $H_{n+1}Q_{*}$.  Since $H_{n+1}Q_{*}$ is FP-injective, there is a map
$F_{n}\mathcolon P_{n}\xrightarrow{}H_{n+1}Q_{*}$ such that
$F_{n}d_{n+1}=\overline{\phi_{n+1}}-E_{n}d_{n+1}$.  Hence
\[
\overline{D_{n}}=E_{n}+F_{n}\mathcolon P_{n} \xrightarrow{}Q_{n+1}/B_{n+1}Q
\]
defines a lift in our commutative square.  

We now choose $D_{n}\mathcolon P_{n}\xrightarrow{}Q_{n+1}$ lifting
$\overline{D_{n}}$, which we can do because $P_{n}$ is projective.  Then
one can easily check that $d_{n+1}D_{n}=\phi_{n}-D_{n-1}d_{n}$, and
also, because $\overline{D_{n}}d_{n+1}=\overline{\phi_{n+1}}$, that
$\phi_{n+1}-D_{n}d_{n+1}$ lands in $B_{n+1}Q$.  This completes the
induction step and the proof.  
\end{proof}

We can now prove Theorem~\ref{thm-class}.

\begin{proof}[Proof of Theorem~\ref{thm-class}]
Suppose the generating hypothesis holds in $\cat{D}(R)$.  In view of
Theorem~\ref{thm-generating}, we need only show that $R$ has weak
dimension at most $1$.  Since $\Tor_{*}(-,M)$ commutes with direct
limits, it suffices to show that the weak dimension of any finitely
presented module is at most $1$.  Since any finitely presented module is
a homology group of a perfect complex, it is enough to show that the
cycles $Z_{n}P$ and the boundaries $B_{n}P$ are flat for all perfect
complexes $P$ and integers $n$.  But $Z_{n}P$ is itself a homology group
of a perfect complex by Corollary~\ref{cor-small}, and so
Theorem~\ref{thm-generating} implies that $Z_{n}P$ is FP-injective.
This means that the short exact sequence
\[
0 \xrightarrow{} Z_{n}P \xrightarrow{} P_{n} \xrightarrow{} B_{n-1}P
\xrightarrow{} 0
\]
is pure.  Now choose a left $R$-module $M$ and apply $-\otimes_{R}M$ to
this short exact sequence.  By purity, it remains exact, and so the
$\Tor$ long exact sequence shows that $\Tor_{1}^{R}(B_{n-1}P,M)=0$.
Since $M$ was arbitrary, $B_{n-1}P$ is flat.  But then $Z_{n}P$, as a
kernel of a surjection of flat modules, is also flat.  

Conversely, assume $R$ has global weak dimension at most $1$ and all
finitely presented $R$-modules are FP-injective.  We need to show that
an arbitrary homology group $M$ of a perfect complex is FP-injective, by
Theorem~\ref{thm-generating}.  By Proposition~\ref{prop-small}, there is
a finitely presented module $F$ and an exact sequence
\[
0 \xrightarrow{} M \xrightarrow{} F \xrightarrow{} F/M \xrightarrow{} 0,
\]
where $F/M$ embeds in a projective module.  Since $R$ has global weak
dimension at most $1$, $F/M$ is flat.  But then the
above exact sequence is pure~\cite[Theorem~4.85]{lam}.  Applying
$\Hom_{R}(N,-)$ to this sequence we get a long exact sequence
\begin{gather*}
0 \xrightarrow{} \Hom_{R}(N,M) \xrightarrow{} \Hom_{R}(N,F)
\xrightarrow{}\Hom_{R} (N, F/M) \\
\xrightarrow{} \Ext^{1}_{R}(N,M)
\xrightarrow{} \Ext^{1}_{R}(N, F) \xrightarrow{} \dotsb 
\end{gather*}
If $N$ is finitely presented, though, the map
$\Hom_{R}(N,F)\xrightarrow{}\Hom_{R}(N,F/M)$ is surjective, since our
original sequence is pure~\cite[Theorem~4.89(5)]{lam}.  By
hypothesis, $\Ext^{1}_{R}(N,F)=0$, so we conclude that
$\Ext^{1}_{R}(N,M)=0$.  Thus $M$ is FP-injective.
\end{proof}

\section{Examples and counterexamples}\label{sec-examples}

In this section, we give conditions under which rings that satisfy the
generating hypothesis must be von Neumann regular, and also give an
example of a ring that satisfies the generating hypothesis yet is not
von Neumann regular, and thus does not satisfy the strong generating
hypothesis.  

\begin{theorem}\label{thm-flat}
A ring $R$ is von Neumann regular if and only if the generating
hypothesis holds in $\cat{D}(R)$ and finitely generated flat submodules
of projective right $R$-modules are projective.
\end{theorem}

\begin{proof}
Assume that the generating hypothesis holds in $\cat{D}(R)$ and finitely
generated flat submodules of projectives are projective.  We will show
that all finitely presented modules, and hence all modules, are flat.
Given a finitely presented module $M$, choose a perfect complex $P$ with
$M\cong H_{n}P$ for some $n$.  We then have a short exact sequence
\[
0 \xrightarrow{} B_{n} P \xrightarrow{} Z_{n} P \xrightarrow{} M
\xrightarrow{} 0.
\]
Now $B_{n}P$ is finitely generated and flat (since it is a submodule of
$P_{n}$) by Theorem~\ref{thm-class}.  By hypothesis, then, $B_{n}P$ is
finitely generated projective.  Hence $B_{n}P$ is FP-injective by
Theorem~\ref{thm-class} again, and so the above exact sequence splits.
Thus $M$ is a summand of $Z_{n}P$, which is flat as well, since it is
also a submodule of $P_{n}$.  So $M$ is flat.  

Conversely, if $R$ is von Neumann regular, then any finitely generated
submodule of a projective module is
projective~\cite[Example~2.32(d)]{lam}.
\end{proof}

This immediately gives the following corollary, implicit
in~\cite{lockridge}.

\begin{corollary}\label{cor-coherent}
A ring $R$ is von Neumann regular if and only if $R$ satisfies the
generating hypothesis and is right coherent.  
\end{corollary}

\begin{proof}
If $R$ is right coherent, then a finitely generated submodule of a
projective module is finitely presented.  If it is also flat, then it is
projective.
\end{proof}

There are a great many rings where finitely generated flat modules are
known to be projective~\cite{puninski-rothmaler}.  The following theorem
contains some cases of this, which are somewhat less satisfactory since
not all von Neumann regular rings satisfy the hypotheses.

\begin{theorem}\label{thm-gen}
Suppose the generating hypothesis holds in $\cat{D}(R)$ and one of the
following hypotheses holds.
\begin{enumerate}
\item $R$ is local \ulp unique maximal right ideal\urp.  
\item $R$ is semiperfect \ulp every finitely generated module has a
projective cover\urp .  
\item $R$ is reduced \ulp no nonzero nilpotents\urp and has finite
uniform dimension \ulp $R$ contains no infinite direct sum of nonzero
right ideals\urp. 
\item $R$ has zero Jacobson radical and finite uniform dimension. 
\item $R$ is right nonsingular \ulp the only element whose right
annihilator is essential in $R$ is $0$\urp and has finite uniform
dimension. 
\item $R$ is simple \ulp no nontrivial two-sided ideals\urp and has
finite uniform dimension.  
\end{enumerate}
Then $R$ is von Neumann regular.  
\end{theorem}

Note that these conditions may not all be independent of each other.  
For example, the authors suspect that if $R$ is both right FP-injective
(as it must be if it satisfies the generating hypothesis) and has finite
uniform dimension, then $R$ may have to be semiperfect.  

\begin{proof}
For a local, semiperfect, or right nonsingular
ring with finite uniform dimension, every finitely generated flat module
is projective; the local case is due to Endo and can be found
in~\cite[Theorem~4.38]{lam}.  The semiperfect case is due to Bass and
is~\cite[Exercise~4.21]{lam}.  The right nonsingular case is due to
Sandomierski~\cite[Corollary~1,p.~228]{sandomierski}.  Every reduced
ring is right nonsingular by~\cite[Lemma~7.8]{lam}; since the
singular elements form a two-sided ideal, every simple ring is also
right nonsingular~\cite[Section~7A]{lam}.  If $R$ is FP-injective,
or in fact only has $\Ext^{1}(R/aR,R)=0$ for all $a\in R$, then having
zero Jacobson radical is equivalent to being right nonsingular,
by~\cite[Theorem~2.1]{nicholson-yousif}.
\end{proof}

Not every von Neumann regular ring has finite uniform dimension.  They
all, however, are right nonsingular~\cite[Corollary~7.7]{lam}.  This
leads to the following theorem. 

\begin{theorem}\label{thm-nonsingular}
A ring $R$ is von Neumann regular if and only if it satisfies the
generating hypothesis, is right nonsingular, and its maximal right ring
of quotients $Q$ is a flat left $R$-module.  
\end{theorem}

The maximal right ring of quotients of $R$ is the endomorphism ring of
the injective hull of $R$ as a right $R$-module, and is much studied in
ring theory.  See~\cite[Section~13]{lam} for an introduction.  When
$R$ is right nonsingular, $Q$ is just equal to the injective hull of
$R$.  

\begin{proof}
Sandomierski~\cite[Theorem~2.9]{sandomierski} proves that if $R$ is
right nonsingular and the maximal right ring of quotients $Q$ is flat as
a left $R$-module, then finitely generated flat submodules of free
$R$-modules (and hence also of projective $R$-modules) are
projective. Theorem~\ref{thm-flat} completes the proof.
\end{proof}

\begin{theorem}\label{thm-counterexample}
There exists a ring $S$ that satisfies the generating hypothesis but is
not von Neumann regular.  
\end{theorem}

Of course, such a ring will not satisfy the strong generating
hypothesis.  Before proving this theorem, we need the following lemmas.

\begin{lemma}\label{lem-weak-dim-cyclic}
Every principal right ideal of a ring $R$ is flat if and only if
whenever $ab=0$ in $R$ there is an $x\in R$ such that $ax=0$ and
$xb=b$.  
\end{lemma}

\begin{proof}
Consider the short exact sequence 
\[
0 \xrightarrow{} \ann_{r} a \xrightarrow{} R \xrightarrow{a \times} aR
\xrightarrow{} 0.
\]
By~\cite[Theorem~4.23]{lam}, $aR$ is flat if and only if for every
$b\in \ann_{r}a$, there is a map $\theta \mathcolon
R\xrightarrow{}\ann_{r}(a)$ with $\theta (b)=b$.  Translating, this
means that $aR$ is flat if and only if whenever $ab=0$, there is an $x$
such that $ax=0$ and $xb=b$.
\end{proof}

\begin{lemma}\label{lem-weak-dim}
A ring $R$ has global weak dimension $\leq 1$ if and only if for every
integer $m$ and every pair of $m\times m$ matrices $A,B$ over $R$ with
$AB=0$, there is an $m\times m$ matrix $X$ over $R$ such that $AX=0$ and
$XB=B$.  
\end{lemma}

\begin{proof}
In view of Lemma~\ref{lem-weak-dim-cyclic}, the matrix condition of this
lemma is equivalent to every principal right ideal of $M_{m}(R)$ being
flat, for all $m\geq 1$. We will use the Morita equivalence between $R$
and $M_{m}(R)$ to prove that this is equivalent to $R$ having global
weak dimension $\leq 1$.  Indeed, if $R$ has global weak dimension $\leq
1$, so does $M_{m}(R)$~\cite[p.~481]{lam}, and so every ideal of
$M_{m}(R)$ is flat.  

Conversely, suppose every principal right ideal of $M_{m}(R)$ is flat
for all $m\geq 1$.  Suppose $I$ is an $m$-generated right ideal of $R$.
Then $I$ corresponds under the Morita equivalence to a principal right
ideal of $M_{m}(R)$~\cite[Remark~17.23(C)]{lam}.  This principal
ideal is flat, and so $I$ is flat as well, since Morita equivalences
preserve flatness~\cite[p.~481]{lam}.

Hence all finitely generated ideals of $R$ are flat.  Since $\Tor$
commutes with direct limits, all ideals of $R$ are flat.  But then $R$
has weak dimension $\leq 1$~\cite[Lemma~4.66]{lam}.  
\end{proof}

\begin{proof}[Proof of Theorem~\ref{thm-counterexample}] We will use the
method of~\cite{prest-rothmaler-ziegler}, who introduce and study
indiscrete rings.  For us, the salient property of indiscrete rings is
that all finitely presented modules over an indiscrete ring are
FP-injective~\cite[Theorem~2.4]{prest-rothmaler-ziegler}.  Thus, we must
find an indiscrete ring that also has weak dimension one.  The
construction given in~\cite[p.~359]{prest-rothmaler-ziegler} begins with
a finite-dimensional algebra $R$ of finite representation type over an
infinite field $F$.  Because we want to end up with something of weak
dimension one, we will take $R$ to have right (and left) global
dimension $1$.  For example, we can take $R$ to be the ring of $2\times
2$ upper triangular matrices over $F$, which is a classical example of a
ring of right (and left) global dimension $1$ that is not von Neumann
regular~\cite[Example~2.36]{lam}.

The method of~\cite{prest-rothmaler-ziegler} is then to construct a map
$\tau \mathcolon R\xrightarrow{}M_{n}R$ and then let $S=R_{\tau }$ be
the direct limit 
\[
S=R_{\tau} = \colim (R \xrightarrow{\tau } M_{n}R\xrightarrow{M_{n}\tau}
M_{n^{2}}R\xrightarrow{M_{n^{2}}\tau} \dotsb )
\]
Then Prest, Rothmaler, and Ziegler show that $S$ is always indiscrete.
Now, in our case, our ring $R$ has global dimension $1$, and therefore
all of the $M_{k}R$ also have global dimension $1$ since they are Morita
equivalent to $R$.  Now, if we take a pair of $m\times m$ matrices $A,B$
over $S$ with $AB=0$, then we can choose $k$ large enough such that
$A,B$ are actually matrices over $M_{n^{k}}R$, and $AB=0$ as such
matrices.  Then Lemma~\ref{lem-weak-dim} shows that there is a matrix
$X$ over $M_{n^{k}}R$, and hence over $S$, with $AX=A$ and $XB=0$.  Thus
Lemma~\ref{lem-weak-dim} implies that $S$ has weak dimension $\leq 1$,
and $S$ cannot be von Neumann regular because $R$ is not
(see~\cite[p.~359]{prest-rothmaler-ziegler}).
\end{proof}

The indiscrete rings of~\cite{prest-rothmaler-ziegler}, of which our
counterexample $S$ is one, have been generalized by Garkusha and
Generalov~\cite{garkusha-generalov} to the class of almost regular
rings, in which all (left or right) finitely presented modules are
FP-injective. The indiscrete rings are the simple almost regular rings.  

We also note that the ring $S$ of Theorem~\ref{thm-counterexample} is in
fact \textbf{weakly semihereditary} in the sense of
Cohn~\cite[p.~13]{cohn}.  This means that if $A$ and $B$ are (not
necessarily square) matrices such that $AB=0$, then there is an
idempotent matrix $E$ such that $AE=A$ and $EB=0$.  Since hereditary
implies weakly semihereditary, each $M_{n}(R)$ in the above proof is
weakly semihereditary, and so the same argument shows that $S$ is as
well.  One can then use (the left module version of)
Lemma~\ref{lem-weak-dim} to see that weakly semihereditary implies
global weak dimension $\leq 1$.

We now turn to some questions we have not been able to answer.  First of
all, the stable homotopy category in topology is more like $\cat{D}(R)$
for a \textbf{graded} ring $R$ (or, better yet, a differential graded
algebra $R$), though, it must be stressed, these categories are still
much simpler than the stable homotopy category.  We have not considered
the generating hypothesis for these $R$.

We could ask whether there is a ring $R$ that satisfies the generating
hypothesis for right $R$-modules but not left $R$-modules.  Such a ring
could not be von Neumann regular, of course.

Also, recall that there is a strongly convergent spectral sequence
whose $E_{2}$ term is $\Ext_{R}^{**}(H_{*}P, H_{*}Q)$ converging to
$\cat{D}(R)(P,Q)_{*}$.  It seems intuitively evident that for the strong
generating hypothesis to hold, this spectral sequence must collapse to
the $0$-line for perfect complexes $P$ and $Q$.  This is in fact true,
since in this case $R$ is von Neumann regular, hence coherent, so the
homology groups $H_{*}P$ are finitely presented modules and therefore
projective.  

However, the situation for $R$ satisfying the generating hypothesis but
not the strong generating hypothesis is less clear.  To satisfy the
generating hypothesis, it must be that every element of
$\Ext^{s,*}(H_{*}P, H_{*}Q)$ with $s>0$ does not survive the spectral
sequence.  But in order not to satisfy the strong generating hypothesis,
there must be an element of $\Hom^{*}(H_{*}P,H_{*}Q)$ for some perfect
$P$ and $Q$ that supports a differential.  It would be intriguing to
understand how this happens.  

Finally, one could define $R$ to satisfy the \textbf{$n$-fold generating
hypothesis} if whenever $f_{1},\dotsc ,f_{n}$ are composable maps of
perfect complexes such that $H_{*}(f_{i})=0$ for all $i$, then
$f_{n}\circ \dotsb \circ f_{1}=0$ in $\cat{D}(R)$.  If we ask for this
condition to hold for all $n$-tuples of composable maps with
$H_{*}f_{i}=0$, not just maps between perfect complexes, then the second
author has shown in his thesis, using work of
Christensen~\cite{christensen}, that $R$ has projective dimension $\leq
n$.  One could then ask for an analogous characterization of rings $R$,
probably in terms of weak dimension, that satisfy the $n$-fold
generating hypothesis, or some strong version of the $n$-fold
generating hypothesis.


\providecommand{\bysame}{\leavevmode\hbox to3em{\hrulefill}\thinspace}
\providecommand{\MR}{\relax\ifhmode\unskip\space\fi MR }
\providecommand{\MRhref}[2]{%
  \href{http://www.ams.org/mathscinet-getitem?mr=#1}{#2}
}
\providecommand{\href}[2]{#2}

\end{document}